\documentstyle[11pt,amssymb,amstex]{amsart}

\def\bn{{\mathbb N}}

\def\l{\lambda} 

\def\m{\mu}

\def\w{\omega} \def\O{\Omega}

\def\nbb{{\mathbf{n}}}
\def\kb{{\mathbf{k}}}

\newtheorem{thm}{Theorem}[section]

\newtheorem{cor}[thm]{Corollary}

\newtheorem{rem}[thm]{Remark}

\theoremstyle{remark}

\def\id{{\bf 1}\!\!{\rm I}}
\def\nb{\nabla}

\begin{document}

\title[Weighted ergodic theorems ]
{Weighted ergodic theorems for  Banach-Kantorovich lattice
$L_{p}(\widehat{\nabla},\widehat{\mu})$}
\author{Inomjon Ganiev}
\address{Inomjon Ganiev\\
 Department of Science in Engineering\\
Faculty of Engineering, International Islamic University Malaysia\\
P.O. Box 10, 50728\\
Kuala-Lumpur, Malaysia}  \email{{\tt inam@@iium.edu.my},{\tt
ganiev1@@rambler.ru}}

\author{Farrukh Mukhamedov}
\address{Farrukh Mukhamedov\\
 Department of Computational \& Theoretical Sciences\\
Faculty of Science, International Islamic University Malaysia\\
P.O. Box, 141, 25710, Kuantan\\
Pahang, Malaysia} \email{{\tt farrukh\_m@@iium.edu.my}, {\tt
far75m@@yandex.ru}}

\begin{abstract}
In the present paper we prove weighted ergodic theorems and
multiparameter weighted ergodic theorems for positive contractions
acting on $L_p(\hat{\nabla},\hat{\mu})$. Our main tool is the use of
methods of measurable bundles of Banach-Kantorovich lattices.
 \vskip 0.3cm \noindent
{\it Mathematics Subject Classification}: 37A30, 47A35, 46B42, 46E30, 46G10.\\
{\it Key words and phrases}: Banach-Kantorovich lattice, positive
contraction, weighted ergodic theorem.

\end{abstract}

\maketitle

\section{Introduction}

The present paper is devoted to the weighted ergodic theorems for
positive contractions acting on Banach-Kantorovich lattice
$L_p(\hat{\nabla},\hat{\mu})$. Note that in \cite{BO,JO} weighted
ergodic theorems for Danford-Schwraz operators acting on
$L_p$-spaces were proved. Further, in \cite{B1,B2} such results were
extended to Banach-valued functions.
 In \cite{LW} it has been considered weighted ergodic
theorems and strong laws of large numbers. In \cite{W} some
properties of the convergence of Banach-valued martingales were
described and their connections with the geometrical properties of
Banach spaces were established too.

It is known that the theory of Banach bundles stemming from the
paper \cite{G2}, where it was proved such a theory has vast
applications in analysis. In  \cite{G1,G2,K1,K2}) the theory of
Banach--Kantorovich spaces were developed. In \cite{Ga3}
Banach-Kantorovich lattice $L_{p}(\widehat{\nabla},\widehat{\mu})$
is represented as a measurable bundle of classical $L_{p}$
--lattices. Hence, with the development of the theory
Banach--Kantorovich spaces
 there naturally arises the necessity to
study some ergodic type theorems for positive contractions and
martingales defined on such spaces. In \cite{CGa} an analog of the
individual ergodic theorem for positive contractions of
Banach-Kantorovich lattices $L_{p}(\widehat{\nabla},\widehat{\mu})$,
has bee established. In \cite{ZC} such a result has been extended to
Orlich-Kantorovich lattices.  In \cite{Ga2} the convergence of
martingales on such lattices is proved. Further, in \cite{GaM} the
"zero-two" law for positive contractions of Banach-Kantorovich
lattice $L_{p}(\widehat{\nabla},\widehat{\mu})$ has been proved.

In the present paper we are going to prove weighted ergodic theorems
and multiparameter weighted ergodic  theorems for positive
contractions acting on $L_p(\hat{\nabla},\hat{\mu})$. We note that
more effective methods to study of Banach-Kantorovich spaces are the
methods of Boolean-valued analysis and measurable bundles (see
\cite{K1},\cite{K2}, \cite{K3}). In the present paper we shall use
the methods of measurable bundles of Banach-Kantorovich lattices.

\section{Preliminaries}

In this section we recall and formulate necessary definitions and
results about Banach-Kantorovich lattices.

Let $(\Omega,\Sigma,\lambda)$ be a measurable space with finite
measure $\lambda$, and $L_0(\Omega)$ be the algebra of all
measurable functions on $\O$ ( here the functions equal a.e. are
identified) and let $\nb(\O)$ be the Boolean algebra of all
idempotents in $L_0(\Omega)$. By $\nb$ we denote an arbitrary
complete Boolean subalgebra of $\nb(\O)$.  By ${\cal
L^{\infty}}(\Omega)$ we denote the set of all measurable essentially
bounded functions on $\O$, and $L^{\infty}(\Omega)$ denote an
algebra of equivalence classes of essentially bounded measurable
functions.

Let $E$ be a linear space over the real field $\mathbb{R}$. By
$\|\cdot\|$ we denote a $L_0(\Omega)$-valued norm on $E$. Then the
pair $(E,\|\cdot\|)$ is called a {\it lattice-normed space (LNS)
over $L_0(\Omega)$}. An LNS $E$ is said to be {\it $d$-decomposable}
if for every $x\in E$ and the decomposition $\|x\|=f+g$ with $f$ and
$g$ disjoint positive elements in $L_0(\Omega)$ there exist $y,z\in
E$ such that $x=y+z$ with $\|y\|=f$, $\|z\|=g$.

Suppose that $(E,\|\cdot\|)$ is an LNS over $L_0(\O)$. A net
$\{x_\alpha\}$ of elements of $E$ is said to be {\it
$(bo)$-converging} to $x\in E$ (in this case we write
$x=(bo)$-$\lim x_\alpha$), if the net $\{\|x_\alpha - x\|\}$
$(o)$-converges to zero in $L_0(\Omega)$ (written as $(o)$-$\lim
\|x_\alpha -x\|=0$). A net $\{x_\alpha\}_{\alpha\in A}$ is called
{\it $(bo)$-fundamental} if $(x_\alpha-x_\beta)_{(\alpha,\beta)\in
A\times A}$ $(bo)$-converges to zero.

An LNS in which every $(bo)$-fundamental net $(bo)$-converges is
called {\it $(bo)$-complete}. A {\it Banach-Kantorovich space (BKS)
over $L_0(\Omega)$} is a $(bo)$-complete $d$-decomposable LNS over
$L_0(\Omega)$. It is well known \cite{K1},\cite{K2} that every BKS
$E$ over $L_0(\Omega)$ admits an $L_0(\Omega)$-module structure such
that $\|fx\|=|f|\cdot\|x\|$ for every $x\in E,\ f\in L_0(\Omega)$,
where $|f|$ is the modulus of a function $f\in L_0(\Omega)$.
 A BKS $({\cal U},\|\cdot\|)$ is called a {\it
Banach-Kantorovich lattice} if  ${\cal U}$ is a vector lattice and
the norm $\|\cdot\|$ is monotone, i.e.  $|u_1|\leq|u_2|$ implies
$\|u_1\|\leq\|u_2\|$. It is known \cite{K1} that the cone ${\cal
U}_+$ of  positive elements is $(bo)$-closed.

Let $(\O,\Sigma,\lambda)$ be the same as above and  $X$ be an
assisting real Banach space $(X(\omega),\|\cdot\|_{X(\omega)})$ to
each point $\omega\in\Omega$, where $ X(\omega)\neq\{0\}$ for all
$\w\in\O$. A {\it section} of $X$ is a function $u$ defined
$\l$-almost everywhere in $\O$ that takes values $u(\w)\in X(\w)$
for all $\w$ in the domain $dom(u)$ of $u$. Let $L$ be a set of
sections. The pair $(X,L)$ is called a {\it measurable Banach bundle
over $\Omega$} if
\begin{enumerate}
   \item[(1)] $\alpha_1 u_1 +\alpha_2 u_2 \in L$ for every $\alpha_1,\alpha_2\in\mathbb{R}$
   and $u_1,u_2\in L$, where\\
 $\alpha_1 u_1 + \alpha_2 u_2 : \omega\in dom(u_1) \cap
dom(u_2) \to \alpha_1 u_1(\omega) + \alpha_2 u_2(\omega)$;
   \item[(2)] the function $\|u\| : \omega\in dom(u) \to
\|u(\omega)\|_{X(\omega)}$ is measurable for every $u\in L$;
    \item[(3)] the set $\{u(\omega) : u\in L, \omega\in dom(u)\}$ is dense in
    $X(\omega)$ for every $\omega\in\O$.
\end{enumerate}

A measurable Banach bundle $(X,L)$ is called {\it measurable bundle
of Banach lattices (MBBL)} if $(X(\omega),\|\cdot\|_{X(\omega)})$ is
a Banach lattice for all $\w\in\O$ and for every $u_1,u_2\in L$ one
has $u_1\vee u_2\in L$, where $u_1 \vee u_2$: $\omega\in{\rm dom}\
(u_1)\cap {\rm dom}\
 (u_2) \rightarrow u_1(\omega)\vee u_2(\omega)$.

A section $s$ is called {\it step-section} if it has a form
$$
s(\omega)=\sum\limits_{i=1}^n \chi_{A_i}(\omega) u_i(\omega),
$$
for some $u_i\in L$, $A_i\in\Sigma$, $A_i\cap A_j=\emptyset$,
$i\neq j$, $i,j=1,\cdots,n$, $n\in\mathbb{N}$, where $\chi_A$ is
the indicator of a set $A$. A section $u$ is called {\it
measurable}  there exists a sequence of step-functions $\{s_n\}$
such that $s_n(\omega)\to u(\w)$ $\l$-a.e.

By $M(\Omega,X)$ we denote the set all measurable sections, and by
$L_0(\Omega,X)$ the factor space of $M(\Omega,X)$ with respect to
the equivalence relation of the equality a.e. Clearly, $L_0(\O,X)$
is an $L_0(\O)$-module. The equivalence class of an element $u\in
M(\O,X)$ is denoted by  $\hat{u}$. The norm of $\hat{u} L_0(\O,X)$
is defined as a class of equivalence in $L_0(\O)$ containing the
function $\|u(\omega)\|_{X(\omega)}$, namely
$\|\hat{u}\|=\widehat{(\|u(\w)\|_{X(\w)})}$. In \cite{G1} it was
proved that $L_0(\Omega,X)$ is a BKS over $L_0(\O)$. Furthermore,
for every BKS $E$ over $L_0(\O)$ there exists a measurable Banach
bundle $(X,L)$ over $\O$ such that $E$ is isomorphic to
$L_0(\Omega)$.

Let $X$ be a MBBL. We put $\hat u\leq\hat v$ if $u(\w)\leq v(\w)$
a.e. One can see that the relation  $\hat u\leq\hat v$ is a partial
order in  $L_0(\Omega,X)$. If $X$ is a MBBL, then  $L_0(\Omega,X)$
is a Banach-Kantarovich lattice \cite{Ga1,Ga3}.

A mapping $\mu : {\nabla}\to L_0(\Omega)$ is called a {\it
$L_0(\Omega)$-valued measure} if the following conditions are
satisfied:
\begin{enumerate}
\item[1)] $\mu(e)\geq 0$ for all $e\in{\nabla}$;
\item[2)] if $\ e\wedge g=0, e,g\in{\nabla},$ then $\mu(e\vee g)=\mu(e)+\mu(g)$;
\item[3)] if $e_n\downarrow 0,\ e_n\in{\nabla},\ n\in\mathbb{N}$, then
$\mu(e_n)\downarrow 0$.
\end{enumerate}

A $L_0(\Omega)$-valued measure $\mu$ is called {\it strictly
positive} if $\allowbreak \mu(e)=\nobreak 0,\ \allowbreak
e\in{\nabla}$ implies $e=0$.

Let a Boolean algebra ${\nabla}(\Omega)$ of all idempotents of
$L_0(\Omega)$ is a regular subalgebra of
 ${\nabla}$.

In the sequel we will consider a strictly positive
$L_0(\Omega)$-valued measure $\m$ with the following property
$\m(ge)=g\m(e)$ for all $e\in\nb$ and $g\in\nb(\O)$.

Let $\nabla_\omega,\ \omega\in\Omega$ be complete Boolean algebras
with strictly positive real-valued measures $\mu_\omega$. Put
 $\rho_\omega(e,g)=\mu_\omega(e\vartriangle g)$,
 $e,g\in\nabla_\omega$. Then $(\nabla_\omega,\mu_\omega)$ is a
 complete metric space. Let us consider a mapping
 $\nabla$, which assigns to each $\omega\in\Omega$ a Boolean algebra
 $\nabla_\omega$. Such a mapping is called a section.

 Assume that $L$ is a nonempty set of sections $\nabla$.
A pair $(\nabla,L)$ is called {\it a measurable bundle of Boolean
algebras over $\O$} if one has
\begin{enumerate}
\item[1)] $(\nabla,L)$  is a measurable bundle of metric spaces (see  \cite{Ga3});
\item[2)] if $e\in L$, then $e^{\perp}\in L$, where $e^{\perp} : \omega\in{\rm dom}\ (e) \rightarrow e^{\perp}(\omega)$;
\item[3)] if $e_1,e_2\in L$, then $e_1\vee e_2\in L$, where
$e_1\vee e_2 : \omega\in{\rm dom}\ (e_1) \cap {\rm dom}\ (e_2)
\rightarrow e_1(\omega) \vee e_2(\omega)$.
\end{enumerate}

Let $M(\Omega,\nabla)$ be the set of all measurable sections, and
$\hat{\nabla}$ be a factorization of $M(\Omega,\nabla)$ with respect
to equivalence relation the equality a.e. Let us define a mapping
$\hat{\mu}:\hat{\nabla}\rightarrow L_0(\Omega)$ by
$\hat{\mu}(\hat{e})=\hat{f}$, where $\hat{f}$ is a class containing
the function $f(\omega)=\mu_\omega(e(\omega))$. It is clear that the
mapping $\hat{\mu}$ is well-defined. It is known that
$(\hat{\nabla},\hat{\mu})$ is a complete Boolean algebra with a
strictly positive $L_0(\Omega)$-valued measure $\hat{\mu}$. Note
that a Boolean algebra ${\nabla}(\Omega)$ of all idempotents of
$L_0(\Omega)$ is identified with a regular subalgebra of
 $\hat{\nabla}$, and one has $\widehat{\mu}(g\widehat{e})=g\widehat{\mu}(\widehat{e})$ for all
 $g\in\nabla(\Omega)$ and $\widehat{e}\in\widehat{\nabla}$.

 The reverse is also tree, namely one has the following

\begin{thm}\cite{Ga3} Let $\tilde{\nabla}$ be a complete Boolean algebra, $\tilde{\mu}$ be a strictly positive
$L_0(\Omega)$-valued measure on $\tilde{\nabla}$, and
$\nabla(\Omega)$ is a regular subalgebra of $\tilde{\nabla}$ and
 $\tilde{\mu}(g\tilde{e})=g\tilde{\mu}(\tilde{e})$ for all $g\in\nabla(\Omega)$, $\tilde{e}\in\tilde{\nabla}$.
 Then there exists a measurable bundle of Boolean algebras
 $(\nabla,L)$ such that $\hat{\nabla}$ is isometrically isomorphic to
 $\tilde{\nabla}$.
 \end{thm}

 By the equality
 $\mu(e)=\frac{\tilde{\mu}(e)}{1+\tilde{\mu}(\id)}$ we define
 $L^\infty(\Omega)$-valued strictly positive measure on
 $\tilde{\nabla}$, where $\id$ is an identity in $\tilde{\nabla}.$

Assume that $p$ is a lifting on $L^\infty(\Omega)$ \cite{G2}.
Define a real-valued quasi-measure on $\tilde{\nabla}$ by
$\mu_\omega^0 (e)=p(\mu(e))(\omega)$ for all $\omega\in\Omega$.

Let $I_\omega^0=\{e\in\tilde{\nabla} : \mu_\omega^0 (e)=0\}$ for all
$\omega\in\Omega$. It is clear that $I_\omega^0$ is an ideal of
$\tilde{\nabla}$. Put $\nabla_\omega^0=\tilde{\nabla} /
 I_\omega^0.$ Then $\nabla_\omega^0$ is a Boolean algebra with strictly positive
 quasi-measure $\mu_\omega^0$. Let us complete the metric space $(\nabla_\omega^0,\rho_\omega)$,
 where $\rho_\omega(e,g)=\mu_\omega^0(e\vartriangle g)$, and completion we denote by
 $\nabla_\omega$. Then $\nabla_\omega$ is a complete Boolean algebra
 with strictly positive real-valued measure
 $\mu_\omega$, which is an extension of $\mu_\omega^0$.

Assume that $\pi_\omega : \tilde{\nabla}\rightarrow\nabla_\omega^0$
is a factor-homomorphism, $i_\omega :
\nabla_\omega^0\rightarrow\nabla_\omega$ is the inclusion, then
$\gamma_\omega = i_\omega\circ\pi_\omega$ is a homomorphism from the
Boolean algebra $\tilde{\nabla}$ into Boolean algebra
$\nabla_\omega$.

Let  $\nabla$ be a mapping, which assigns to each $\omega\in\Omega$
a Boolean algebra $\nabla_\omega$ with strictly positive real-valued
measure $\mu_\omega$, and $L=\{\tilde{e} :
 \tilde{e}(\omega)=\gamma_\omega(e), e\in\tilde{\nabla}\}$.
It is known \cite{Ga3} that the pair $(\nabla,L)$ is a measurable
bundle of Boolean algebras and $\hat{\nabla}$ is isometrically
isomorphic to $\tilde{\nabla}$. Moreover, one has
$\widehat{\mu}=\tilde{\mu}$ (see \cite{Ga3} for more details).

By $L_0(\hat{\nabla},\hat{\mu})$ we denote an order complete vector
lattice $C_{\infty}(Q(\hat{\nabla})),$ where $Q(\hat{\nabla})$ is
the Stonian compact associated with complete Boolean algebra
$\hat{\nabla}$.  For $\widehat{f},\widehat{g}\in
L_0(\hat{\nabla},\hat{\mu})$ we let
 $\widehat{\rho}(\widehat{f},\widehat{g})=\allowbreak
 \int \frac{|\hat{f}-\hat{g}|}{1+|\hat{f}-\hat{g}|}d\hat{\mu}.$
 Then it is know \cite{Ga3} that $\widehat{\rho}$ is an
 $L_0(\Omega)$-valued metric on  $L_0(\hat{\nabla},\hat{\mu})$ and
 $(L_0(\hat{\nabla},\hat{\mu}),\widehat{\rho})$ is isometrically isomorphic to the measurable bundle
 of metric spaces $L_0(\nabla_\omega,\mu_\omega),$ where $\rho_\omega(a,b)= \int
 \frac{|a-b|}{1_{\omega}+|a-b|}d\mu_\omega.$ In particularly, each element
 $\widehat{f}\in L_0(\hat{\nabla},\hat{\mu})$ can be identified with the measurable section $\{f(\omega)\}_{\omega\in
 \Omega},$ here $f(\omega)\in L_0(\nabla_\omega,\mu_\omega).$

Following  the well known scheme of the construction of
$L_p$-spaces,  a space $L_p(\nb,\m)$ can be defined by
$$
L_p(\hat{\nb},\hat{\m})=\left\{\hat f\in L_0(\hat\nabla):
\int|\hat f|^pd\hat{\m} - \textrm{exist} \ \right\}, \ \ \ p\geq 1
$$
where $\hat\m$ is a  $L_0(\O)$-valued measure on $\hat\nb$.

It is known \cite{K1} that $L_p(\hat\nb,\hat\m)$ is a BKS over
$L_0(\O)$ with respect to the $L_0(\O)$-valued norm $\|\hat
f\|_{L_p(\hat{\nb},\hat{\m})}=\bigg(\int|\hat
f|^pd\hat{\m}\bigg)^{1/p}$. Moreover, $L_p(\hat{\nb},\hat{\m})$ is a
Banach-Kantorovich lattice (see \cite{K2},\cite{Ga3}).

Let $X$ be a mapping assisting an $L_p$-space constructed by a
real-valued measure $\m_\omega$, i.e.
$L_p(\nabla_\omega,\mu_\omega)$ to each point $\omega\in\O$ and let
$$
L=\bigg\{\sum\limits_{i=1}^n \alpha_i e_i : \alpha_i\in
{\mathbb{R}}, \ \ e_i\in M(\Omega,\nabla),\ i=\overline{1,n},\
n\in\mathbb{N}\bigg\}$$ be a set of sections. In \cite{Ga3, GaC}
it has been established that the pair $(X,L)$ is a measurable
bundle of Banach lattices and $L_0(\Omega,X)$ is modulo ordered
isomorphic to $L_p(\hat\nabla,\hat\mu)$.

Let as before $p\geq 1$ and $L_p(\hat{\nabla},\hat{\mu})$ be a
Banach-Kantorovich lattice, and $L_p(\nabla_\omega,\mu_\omega)$ be
the corresponding $L_p$-spaces constructed by a real valued
measures. Let $T:L_p(\hat{\nabla},\hat{\mu})\to
L_p(\hat{\nabla},\hat{\mu})$ be a linear mapping. As usual  we will
say that $T$ is {\it positive} if $T\hat{f}\geq 0$ whenever
$\hat{f}\geq 0$. We say that $T$ is a {\it $L_0(\Omega)$-bounded
mapping} if there exists a function $k\in L_0(\Omega)$ such that
$\|T\hat{f}\|_{L_p(\hat\nb,\hat\m)}\leq k
\|\hat{f}\|_{L_p(\hat\nb,\hat\m)}$ for all $\hat{f}\in
L_p(\hat\nabla,\hat\mu)$. For a such mapping we can define an
element of $L_0(\O)$ as follows
$$
\|T\| =\sup\limits_{\|\hat{f}\|_{L_p(\hat\nb,\hat\m)}\leq\id}
\|T\hat{f}\|_{L_p(\hat\nb,\hat\m)},
$$
which is called an {\it $L_0(\O)$-valued norm} of $T$. If
$\|T\hat{f}\|_{L_p(\hat\nb,\hat\m)}\leq
\|\hat{f}\|_{L_p(\hat\nb,\hat\m)}$ then a mapping $T$ is said to
be a {\it $L_p(\hat\nb,\hat\m)$ contraction}.

The set of all essentially bounded functions w.r.t. $\hat{f}$
taken from $L_0(\hat\nb,\hat\m)$ is denoted by
$L^\infty(\hat\nb,\hat\m)$.

Let $\hat{\nabla}^{(1)}$ be a regular Boolean subalgebra of
$\hat{\nabla}$, and $\hat{\mu}^1$ is the restriction of
$\hat{\mu}$ onto $\hat{\nabla}^1$. Then according to Theorem 4.2.9
\cite{K1} there exists the {\it conditional expectation} operator
$E(\cdot | \hat{\nabla}^1) :
L_1(\hat{\nabla},\hat{\mu})\rightarrow
L_1(\hat{\nabla}^{(1)},\hat{\mu}^1)$  which satisfies the
following conditions:
\begin{enumerate}
\item[(a)] $E(\cdot | \hat{\nabla}^{(1)})$ is linear, positive and idempotent;

\item[(b)] for every  $\hat{f}\in
 L_1(\hat{\nabla},\hat{\mu})$ one has $\int E(\hat{f} | \hat{\nabla}^{(1)}) d\hat{\mu} = \int \hat{f} d\hat{\mu}$;

\item[(c)] $E(\hat{g}\hat{f} | \hat{\nabla}^{(1)}) = \hat{g}
E(\hat{f}
 | \hat{\nabla}^{(1)})$ for every $\hat{g}\in
 L^\infty(\hat{\nabla}^{(1)},\hat{\mu}^{(1)})$ and $\hat{f}\in
 L_1(\hat{\nabla},\hat{\mu})$.

\item[(d)] Moreover, one has $\|E(\hat{f} | \hat{\nabla}^{(1)})
\|_{L_{1}(\widehat{\nabla},\widehat{\mu})} \leq
 \|\hat{f}\|_{L_{1}(\widehat{\nabla},\widehat{\mu})}$
 for every $\hat{f}\in L_1(\hat{\nabla},\hat{\mu})$ and $E(\id| \hat{\nabla}^{(1)}) = \id$.
\end{enumerate}

In the sequel we will need the following

\begin{thm}\label{T-o}\cite{GaM,Ga3} Let $T : L_p(\hat{\nabla},\hat{\mu})\to
L_p(\hat{\nabla},\hat{\mu})$ be a positive linear
$L_p(\hat{\nabla},\hat{\mu})$ contraction such that $T\id\leq\id$.
Then for every $\omega\in\Omega$ there exists a positive
$L_p(\nabla_\omega,\mu_\omega)$ contraction $T(\omega ):
L_p(\nabla_\omega,\mu_\omega)\to L_p(\nabla_\omega,\mu_\omega)$
such that $T(\omega) f(\omega) = (T\hat{f})(\omega)$ $\l$-a.e. for
every $\hat{f}\in L_p(\hat{\nabla},\hat{\mu})$.
\end{thm}

By means of measurable bundle of $L_p$-spaces and at each bundle
applying classical ergodic theorem, it has been proved the following

\begin{thm}\label{erg-t}\cite{CGa} Let $p>1,\ q>1$ with $\frac{1}{p}+\frac{1}{q}=1$ and
$T : L_p(\hat{\nabla},\hat{\mu})\rightarrow
L_p(\hat{\nabla},\hat{\mu})$ be a linear positive
$L_p(\hat{\nabla},\hat{\mu})$ contraction with $T\id\leq\id$. Then
for every $\hat{f}\in L_p(\hat{\nabla},\hat{\mu})$ one has
\begin{enumerate}
\item[(i)] the sequence
$$
s_n(|\hat{f}|)=\frac{1}{n}
 \sum\limits_{i=0}^{n-1} T^i (|\hat{f}|)
$$
is bounded in $L_p(\hat{\nabla},\hat{\mu})$, and
$$\|\sup\limits_{n\geq 1} s_n(|\hat{f}|)\|_{L_{p}(\widehat{\nabla},\widehat{\mu})} \leq q
\|\hat{f}\|_{L_{p}(\widehat{\nabla},\widehat{\mu})};$$

\item[(ii)] there exists an element $\tilde{f}\in
 L_p(\hat{\nabla},\hat{\mu})$ such that the sequence $s_n(\hat{f})$ $(o)$-converges to $\tilde{f}$
 in $L_p(\hat{\nabla},\hat{\mu})$.
\end{enumerate}
\end{thm}

It is worth to mention that in \cite{ZC} it has been proved that for
every $\hat{f}\in L_p(\hat{\nabla},\hat{\mu})$ at $p\geq1$ the
averages $s_n(\hat{f})$ $(bo)$-converge in
$L_p(\hat{\nabla},\hat{\mu})$.

\section{$(o)$-convergence}

In this section we provide some auxiliary facts related to
$(o)$-convergence of sequence $\hat{f}_\nbb$ from
$L_0(\hat{\nabla},\hat{\mu})$ and $(o)$-convergence of the
sequence $\{f_\nbb(\omega)\}$ from
$L_0(\nabla_\omega,\mu_\omega)$.

For $\nbb=(n_1,\dots,n_d)\in\bn^d$ denote
$m(\nbb)=\min\{n_1,\dots,n_d\}$. In the sequel by $\nbb\to\infty$ we
mean $m(\nbb)\to\infty$. For $\nbb=(n_1, n_2,\dots, n_d), \kb=(k_1,
k_2,\dots, k_d)$ we write ${\mathbf{n}}\leq{\mathbf{k}}$ if and only
if $n_i\leq k_i$ (i=1,2,\dots,d), $\nbb<\mathbf{k}$ if $n_i< k_i$
(i=1,2,\dots,d).

\begin{thm}\label{2.3} Let $\widehat{f}_{\mathbf{n}}\in
L_0(\hat{\nabla},\hat{\mu}).$ Then $\sup\limits_\nbb
 \hat{f}_\nbb$ exists in $L_0(\hat{\nabla},\hat{\mu})$ if and only
 if $\sup\limits_\nbb
 {f}_\nbb(\omega)$ exists in $L_0(\nabla_\omega,\mu_\omega)$ for
 a.e. $\omega\in\Omega$. In the later case, one has
 $(\sup\limits_\nbb
 \hat{f}_\nbb)(\omega)=\sup\limits_\nbb
 {f}_\nbb(\omega)$ for
 a.e. $\omega\in\Omega$.
\end{thm}

\begin{pf} Assume that  $g(\omega)=\sup\limits_\nbb
 {f}_\nbb(\omega)$ exists in $L_0(\nabla_\omega,\mu_\omega)$ for
 a.e. $\omega\in\Omega$. Denote
 $\hat{g}_\nbb=\sup\limits_{{\mathbf{1}}
 \leq{\mathbf{k}}\leq{\mathbf{n}}}\hat{f}_\kb$ in $L_0(\hat{\nabla},\hat{\mu}).$
 Then ${g}_\nbb(\omega)=\sup\limits_{{\mathbf{1}}
 \leq{\mathbf{k}}\leq{\mathbf{n}}}{f}_\kb(\omega)$ for
 a.e. $\omega\in\Omega$.

 Obviously, that ${g}_{\mathbf{n}}(\omega)\uparrow{g}(\omega)$ as
 ${\mathbf{n}}\rightarrow\infty$ for
 a.e. $\omega\in\Omega$. The relation ${g}_{\mathbf{n}}(\omega)\uparrow{g}(\omega)$ implies that
 ${g}_{\mathbf{n}}(\omega)\stackrel{\rho_\omega}{\rightarrow}
 {g}(\omega)$ for
 a.e. $\omega\in\Omega$, this means $g\in  M(\Omega,X)$ and
 $\hat{g}\in L_0(\hat{\nabla},\hat{\mu})$.

Let us prove that $\hat{g} =
 \sup\limits_{\mathbf{n}} \hat{f}_{\mathbf{n}}$ in  $L_0(\hat{\nabla},\hat{\mu})$. It is
 clear that ${g}(\omega)\geq {f}_{\mathbf{n}}(\omega)$ for
 a.e. $\omega\in\Omega$. Therefore, $\widehat{{g}}\geq
 \widehat{{f}}_{\mathbf{n}}$ for all $\nbb\in
\bn^d$.

Let $\widehat{\varphi}\in L^0(\hat{\nabla},\hat{\mu})$ and
$\widehat{{\varphi}}\geq
 \widehat{{f}}_{\mathbf{n}}$ for all ${\mathbf{n}}\in
\bn^d.$ Then ${\varphi}(\omega)\geq {f}_{\mathbf{n}}(\omega)$ for
any ${\mathbf{n}}\in \bn^d$. Hence, ${\varphi}(\omega)\geq
{g}(\omega),$ for a.e. $\omega\in\Omega$, i.e. $\hat{\varphi}\geq
 \hat{g}$. This yields that $\hat{g} = \sup\limits_{{\mathbf{n}}\in\bn^d}
 \hat{f}_{\mathbf{n}}$.

Conversely, let us assume that there exists such $\hat{\psi}\in
 L_0(\hat{\nabla},\hat{\mu})$ that $\hat{\psi} = \sup\limits_{\mathbf{n}\geq
 1} \hat{f}_{\mathbf{n}} = \sup\limits_{{\mathbf{n}}\in\bn^d}
 \hat{g}_{\mathbf{n}}$.

From  $g_{\mathbf{n}}(\omega) =
 \sup\limits_{1\leq {\mathbf{k}} \leq {\mathbf{n}}} f_{\mathbf{k}}(\omega)$ for  a.e.
 $\omega\in\Omega$, we find $\psi(\omega)\geq f_{\mathbf{n}}(\omega)$ for all ${\mathbf{n}}\in\bn^d$
 for  a.e. $\omega\in\Omega$.
 Hence, one gets $\psi(\omega)\geq \sup\limits_{{\mathbf{n}}\in\bn^d} f_{\mathbf{n}}(\omega)
 = \sup\limits_{{\mathbf{n}}\in\bn^d} g_{\mathbf{n}}(\omega)$ for  a.e. $\omega\in\Omega$ .
 As $\hat{g}_{\mathbf{n}} \rightarrow\hat{\psi}$ in
 metric $\hat{\rho}$, then  $g_{\mathbf{n}}(\omega)\rightarrow\psi(\omega)$ in
 metric $\rho_\omega$ for  a.e.
 $\omega\in\Omega$.
 Since $\{g_{\mathbf{n}}(\omega)\}$ is increasing then $\psi(\omega) =
 \sup\limits_{{\mathbf{n}}\in\bn^d} g_{\mathbf{n}}(\omega)$ for  a.e.
 $\omega\in\Omega$.
\end{pf}

 From this theorem immediately follows two corollaries.

\begin{cor}\label{2.4} Let
 $\{\hat{f}_{\mathbf{n}}\}\subset L_0(\hat{\nabla},\hat{\mu})$. Then
$\inf\limits_{{\mathbf{n}}\in\bn^d} \hat{f}_{\mathbf{n}}$ exists in
$L_0(\hat{\nabla},\hat{\mu})$
 if and only
 if $\inf\limits_{{\mathbf{n}}\in\bn^d}
 f(\omega)$ exists in $L_0(\nabla_\omega,\mu_\omega)$ for  a.e.
 $\omega\in\Omega$. In later case, one has $(\inf\limits_{{\mathbf{n}}\in\bn^d}
 \hat{f}_{\mathbf{n}})(\omega)= \inf\limits_{{\mathbf{n}}\in\bn^d} f_{\mathbf{n}}(\omega)$ for  a.e.
 $\omega\in\Omega$.
 \end{cor}

\begin{cor}\label{2.5} Let $\hat{f}_{\mathbf{n}} \in
 L_0(\hat{\nabla},\hat{\mu})$. If $\hat{f}_{\mathbf{n}} \stackrel{(o)}\rightarrow \hat{f}$
 for some $\hat{f}\in L_0(\hat{\nabla},\hat{\mu})$, then
 $f_{\mathbf{n}}(\omega) \stackrel{(o)}\rightarrow f(\omega)$ in $L_0(\nabla_\omega,\mu_\omega)$
 for  a.e. $\omega\in\Omega$. Conversely, if $f_{\mathbf{n}}(\omega)
  \stackrel{(o)}\rightarrow g(\omega)$ for some $g(\omega)\in
 L_0(\nabla_\omega,\mu_\omega)$ for  a.e. $\omega\in\Omega$,
 then $\hat{g}\in L_0(\hat{\nabla},\hat{\mu})$ and $\hat{f}_{\mathbf{n}} \stackrel{(o)}\rightarrow
 \hat{g}$ in $L_0(\hat{\nabla},\hat{\mu})$.\\
\end{cor}

\section{Weighted ergodic theorems}

In this section we shall prove some weighted ergodic theorems in
$L_p(\hat{\nabla},\hat{\mu})$.

First we recall that a sequence $\{\alpha(k)\}$ is called {\it
Besicovich} if for every $\varepsilon>0$ there is a sequence of
trigonometric polynomials $\psi_{\varepsilon}$, such that
$$\lim\limits_{N\rightarrow \infty}\sup\frac{1}{N}\sum\limits_{k=1}^{N-1}|\alpha(k)-\psi_{\varepsilon}(k)|<\varepsilon$$

We say that $\{\alpha(k)\}$ is {\it bounded Besicovich} if
$\alpha(k)\in \ell^{\infty}.$ In what follows, we consider only
bounded, real Besicovich sequences.

\begin{thm}\label{B-erg1} Let $\ T :
L_1(\hat{\nabla},\hat{\mu})\rightarrow
L_1(\hat{\nabla},\hat{\mu})$ be a positive linear
$L_1(\hat{\nabla},\hat{\mu})$ contraction with $T\id\leq\id$, and
$\{\alpha(k)\}$ be a bounded Besicovich sequence. Then for every
$\hat{f}\in L_1(\hat{\nabla},\hat{\mu})$ the averages
$$\widetilde{A_{N}}(\widehat{f})=\frac{1}{N}\sum\limits_{k=1}^{N-1}\alpha(k)T^k\widehat{f}$$
$(o)$-converge in $L_0(\widehat{\nabla},\widehat{\mu}).$
\end{thm}

\begin{pf}
According to Theorem \ref{T-o} for each $\omega\in\Omega$ there
exists a positive contraction $T_\omega :
 L_1(\nabla_\omega,\mu_\omega)\rightarrow L_1(\nabla_\omega,\mu_\omega)$,
such that $T_\omega f(\omega) = (T\hat{f})(\omega)$ for every
$\hat{f}\in L_1(\hat{\nabla},\hat{\mu})$ and a.e. $\omega\in\Omega$.
Then for every $\hat{f}\in L_1(\hat{\nabla},\hat{\mu})$ we have
\begin{eqnarray}\label{A_n}
\widetilde{A_{N}}(\widehat{f})(\omega)&=&(\frac{1}{N}\sum\limits_{k=1}^{N-1}\alpha(k)T^k\widehat{f})(\omega)\nonumber\\
&=&\frac{1}{N}\sum\limits_{k=1}^{N-1}\alpha(k)T^k(\omega){f}(\omega)={A_{N}}(\widehat{f}(\omega))
\end{eqnarray}
for a.e. $\omega\in \Omega$.

By means of Theorem 1.4 \cite{JO} one gets the existence of the
limit
$$(o)-\lim\frac{1}{N}\sum\limits_{k=1}^{N-1}\alpha(k)T^k(\omega){f}(\omega)$$
for every $\hat{f}\in L_1(\hat{\nabla},\hat{\mu})$, and the limit
belongs to $L_0(\nabla_\omega,\mu_\omega)$ for a.e. $\omega\in
 \Omega$. This means
\begin{equation*}
\widetilde{A_{N}}(\widehat{f})(\omega)={A_{N}}(\widehat{f}(\omega)){\stackrel{(o)}{\rightarrow}}\widetilde{f}(\omega)
\end{equation*}
in $ L_0(\nabla_\omega,\mu_\omega)$ for a.e. $\omega\in \Omega$ and
for some $\widetilde{f}(\omega)\in L_0(\nabla_\omega,\mu_\omega).$
Due to Corollary \ref{2.5} we obtain
\begin{equation*}
\widetilde{A_{N}}(\widehat{f})=\widehat{{A_{N}}(\widehat{f}(\omega)}){\stackrel{(o)}{\rightarrow}}\widehat{f}
=\widehat{\widetilde{f}(\omega)}
\end{equation*}
in $L_0(\widehat{\nabla},\widehat{\mu}).$
\end{pf}

\begin{rem} In case $\alpha(k)=1$ Theorem \ref{B-erg1} implies Theorem 3.2 (i) \cite{ZC}
\end{rem}

\begin{cor} Let $\ T :
L_1(\hat{\nabla},\hat{\mu})\rightarrow
L_1(\hat{\nabla},\hat{\mu})$ be a positive linear
$L_1(\hat{\nabla},\hat{\mu})$ contraction with $T\id\leq\id$, and
$j_k$ be an increasing sequence of positive integers such that
$\sup\limits_k\frac{j_k}{k}<\infty.$ Then for every $\hat{f}\in
L_1(\hat{\nabla},\hat{\mu})$ the limit
$$(o)-\lim\frac{1}{N}\sum\limits_{k=1}^{N-1}T^{j_k}\widehat{f}$$
exists in  $L_0(\widehat{\nabla},\widehat{\mu})$.
\end{cor}

\begin{pf} Let us define $$\alpha(k)=\left\{%
\begin{array}{ll}
    0, & \hbox{ if $k\neq j_k$;} \\
    1, & \hbox{ if $k=j_k$.} \\
\end{array}%
\right.$$ Then is known \cite{R} that $\alpha(k)$ is a Besicovich
sequence. Hence, Theorem \ref{B-erg1} implies the assertion.
\end{pf}

\begin{thm}\label{B-erg2} Let $p>1,\ q>1$ with  $\frac{1}{p}+\frac{1}{q}=1$, and
$T : L_p(\hat{\nabla},\hat{\mu})\rightarrow
L_p(\hat{\nabla},\hat{\mu})$ be a positive linear
$L_1(\hat{\nabla},\hat{\mu})$ contraction with $T\id\leq\id$, and
$\{\alpha(k)\}$ be a bounded Besicovich sequence.
 Then for every $\hat{f}\in L_p(\hat{\nabla},\hat{\mu})$ one has
\begin{enumerate}
\item[(i)] the sequence $A_N(|\hat{f}|)$ is bounded in
$L_p(\hat{\nabla},\hat{\mu})$ and
$$\|\sup\limits_{N\geq 1} A_N(|\hat{f}|)\|_{L_{p}(\widehat{\nabla},\widehat{\mu})} \leq
q \sup\limits_{k}|\alpha(k)|\|\hat{f}\|_{L_{p}(\widehat{\nabla},\widehat{\mu})};
$$

\item[(ii)] there exists an element $\tilde{f}\in
 L_p(\hat{\nabla},\hat{\mu})$, such that the sequence  $A_N(\hat{f})$ $(o)$-converges to $\tilde{f}$
 in $L_p(\hat{\nabla},\hat{\mu})$.
\end{enumerate}
\end{thm}

\begin{pf} (i) Due to Theorem \ref{T-o} for each $\hat{f}\in
L_p(\hat{\nabla},\hat{\mu})$ and for a.e.  $\omega\in \Omega$ one
has $\widetilde{A_{N}}(\widehat{f})(\omega)={A_{N}}({f}(\omega))$
(see \eqref{A_n}).  Applying  Akcoglu's Theorem \cite{Kr} we have,
that $\sup\limits_{N\geq 1} A_N(|{f(\omega)}|)\in
L_p(\nabla_\omega,\mu_\omega)$ and
$$\|\sup\limits_{N\geq 1}
A_N(|{f(\omega)}|)\|_{L_p(\nabla_\omega,\mu_\omega))} \leq q
\sup_{k}|\alpha(k)|\|{f(\omega)}\|_{L_p(\nabla_\omega,\mu_\omega)}$$
for a.e. $\omega\in \Omega.$  Then we get
\begin{eqnarray*}
\|\sup\limits_{N\geq 1}
A_N(|\hat{f}|)\|_{L_{p}(\widehat{\nabla},\widehat{\mu})}&=&\widehat{{\|\sup\limits_{N\geq
1} A_N(|{f(\omega)}|)\|}_{L_p(\nabla_\omega,\mu_\omega))} }\\
&\leq&
q\sup\limits_{k}|\alpha(k)|\|\widehat{{f(\omega)}}\|_{L_p(\nabla_\omega,\mu_\omega)}\\
&=&q\sup\limits_{k}|\alpha(k)|{\|\widehat{f}\|}_{L_{p}(\widehat{\nabla},\widehat{\mu})}.
\end{eqnarray*}

(ii) Since, $T$ is positive linear $L_1(\hat{\nabla},\hat{\mu})$
contraction in $L_p(\nabla_\omega,\mu_\omega)$ and $T\id\leq\id$,
then by the Theorem 2.2 $T(\omega)$ is
$L_1(\nabla_\omega,\mu_\omega))$ contraction in
$L_p(\nabla_\omega,\mu_\omega))$ and
$T(\omega)\id(\omega)\leq\id(\omega)$. This means that $T(\omega)$
is
$L_1(\nabla_\omega,\mu_\omega))-L_\infty(\nabla_\omega,\mu_\omega))$
contraction $L_p(\nabla_\omega,\mu_\omega))$ for a.e. $\omega\in
\Omega.$ From Theorem 1.2 \cite{B1} we find that the sequence
${A_{N}}({f}(\omega))$ $(o)$- converges to some limit
$\widetilde{f}(\omega)$ a.e. $\omega\in \Omega,$ for every
${f(\omega)}\in L_p(\nabla_\omega,\mu_\omega)$. Then Corollary
\ref{2.5} implies that
$$\widetilde{A_{N}}(\widehat{f})=\widehat{A_{N}({f}(\omega))}{\stackrel{(o)}{\rightarrow}} \widehat{\widetilde{f}(\omega)}
=\widetilde{f}$$ in $L_0(\widehat{\nabla},\widehat{\mu}).$  Due to
$$\|\sup\limits_{N\geq 1}
A_N(|\hat{f}|)\|_{L_{p}(\widehat{\nabla},\widehat{\mu})}\leq q
\sup\limits_{k}|\alpha(k)|{\|\widehat{{f}}\|}_{L_{p}(\widehat{\nabla},\widehat{\mu})}$$
one finds
$$\sup\limits_{N\geq 1} A_N(|\hat{f}|)\in
L_{p}(\widehat{\nabla},\widehat{\mu}).$$ Therefore,
$$A_N(\hat{f}){\stackrel{(o)}{\rightarrow}}\widetilde{f}$$ in
$L_p(\widehat{\nabla},\widehat{\mu}).$
\end{pf}

\section{Multiparameter weighted ergodic theorems}

In what follows, given $\nbb=(n_1,n_2,\dots,n_d)$ we denote
$|{\mathbf{n}}|=n_1\cdot n_2\cdots n_d$, and
${\mathbf{1}}=(1,1,\dots,1).$

Let $T_1,T_2,\dots,T_d$ be $d$ linear positive
$L_p(\hat{\nabla},\hat{\mu})$ contractions in
$L_p(\hat{\nabla},\hat{\mu})$, $1<p<\infty$. Then we denote
${\mathbf{T}}^\nbb=T_1^{n_1}\cdots T_d^{n_d}$, where
$\nbb=(n_1,n_2,\dots,n_d)\in\bn^d$.

\begin{thm}\label{Ind1} Let ${\mathbf{T}}=(T_1,T_2,\dots,T_d)$
denoted $d$  linear  positive $L_p(\hat{\nabla},\hat{\mu})$
contractions in $L_p(\hat{\nabla},\hat{\mu})$, $1<p<\infty$ such
that $T_i\id\leq\id$ for all $i:$ $1\leq i\leq d.$ Then for every
$\widehat{f}\in L_p(\hat{\nabla},\hat{\mu})$ one has
\begin{enumerate}
\item[(i)] The averages
$$
S_\nbb(|\widehat{f}|)=\frac{1}{|\nbb|}
 \sum\limits_{\kb=1}^{\nbb}{\mathbf{T}}^\kb(|\widehat{f}|)
$$
 is bounded in $L_p(\hat{\nabla},\hat{\mu})$, and one has
 $$
 {\|\sup\limits_{{\mathbf{n}}}S_\nbb(|\widehat{f}|)\|}_{L_{p}(\widehat{\nabla},\widehat{\mu})}
 \leq q^d {\|\widehat{f}\|}_{L_{p}(\widehat{\nabla},\widehat{\mu})};$$

\item[(ii)] There exists an element $\widetilde{f}\in
 L_p(\hat{\nabla},\hat{\mu}),$ such that $S_\nbb(\widehat{f})$
 $(o)$- convergence to $\widetilde{f}$ in $L_p(\hat{\nabla},\hat{\mu}).$
\end{enumerate}
\end{thm}

\begin{pf} Due to Theorem \ref{T-o} one has
\begin{eqnarray*}
S_\nbb(\widehat{f})(\omega)&=&\frac{1}{|\nbb|}
 \sum\limits_{k_1=1}^{n_1}\dots\sum\limits_{k_d=1}^{n_d}T_1^{k_1}\cdots T_d^{k_d}(\widehat{f}))(\omega)\\[2mm]
&=& \frac{1}{|\nbb|}
\sum\limits_{k_1=1}^{n_1}\cdots\sum\limits_{k_d=1}^{n_d}T_1^{k_1}(\w)\cdots
T_d^{k_d}(\omega)
 (\widehat{f}(\omega))\\[2mm]
 &=&S_\nbb(\omega)(f(\omega))
\end{eqnarray*}
for  a.e. $\omega\in\Omega.$

By Theorem 1.2 \cite{Kr} (p.196) we have
$\|\sup\limits_{\nbb}S_\nbb(\omega)(|f(\omega)|)\|_{L_{p}(\widehat{\nabla},\widehat{\mu})}\leq
q^d \|f(\omega)\|_{L_{p}(\widehat{\nabla},\widehat{\mu})}$ and
there is $g(\omega)\in L_p(\nabla_\omega,\mu_\omega)$ such that
 $S_{\mathbf{n}}(\omega)(f(\omega))\stackrel{(o)}\rightarrow
 g(\omega)$ in $L_p(\nabla_\omega,\mu_\omega)$. According Theorem \ref{2.3}
  $\sup\limits_{{\mathbf{n}}}S_\nbb(\widehat{f})$ exists in $L_p(\hat{\nabla},\hat{\mu})$,
  and one has
 $${\|\sup\limits_{{\mathbf{n}}}S_\nbb(|\widehat{f}|)\|}_{L_{p}(\widehat{\nabla},\widehat{\mu})}
 \leq q^d {\|\widehat{f}\|}_{L_{p}(\widehat{\nabla},\widehat{\mu})}.$$

 By Corollary \ref{2.5} we obtain that
$$S_\nbb(\widehat{f})=\widehat{S_\nbb(\omega)(f(\omega))}\stackrel{(o)}\rightarrow
 \widehat{g(\omega)}=\widetilde{f}$$ in
 $L_0(\hat{\nabla},\hat{\mu})$. Since $S_\nbb(\widehat{f})$
 is bounded in $L^p(\hat{\nabla},\hat{\mu})$, then
 $S_\nbb(\widehat{f})\stackrel{(o)}\rightarrow\widetilde{f}$
 in $L_p(\hat{\nabla},\hat{\mu})$.
\end{pf}

The next theorem is an analog of multiparameter weighted individual
ergodic theorem in Banach--Kantorovich lattice
$L_p(\hat{\nabla},\hat{\mu}).$

Recall that (see \cite{JO}) a class of weights $\{\alpha(\kb):\
\kb\in \bn^d\}$ is called {\it Besicovich}, if for any
$\varepsilon>0$ there is sequence of trigonometric polynomials in
$d$ variables , $\psi_\varepsilon$ such that
$$\limsup\limits_{\mathbf{n}\rightarrow\infty}\frac{1}{\mathbf{|n|}}
\sum\limits_{{\mathbf{k}}=1}^{{\mathbf{n}}}|\alpha(\mathbf{k})-\psi_\varepsilon(\mathbf{k})|<\varepsilon.$$

We say that $\{\alpha(\mathbf{k})\}$ is {\it bounded Besicovich}
if $\alpha(\mathbf{k})\in \ell^{\infty}.$

\begin{thm}\label{ergB3} Let ${\mathbf{T}}=(T_1,T_2,\dots,T_d)$
denoted $d$ linear  positive $L_1(\hat{\nabla},\hat{\mu})$
contractions in $L_p(\hat{\nabla},\hat{\mu})$, $1<p<\infty$ such
that $T_i\id\leq\id$ for all $i:$ $1\leq i\leq d$, and
$\alpha(\mathbf{k})$ be a bounded Besicovich weights.
  Then for every $\widehat{f}\in L_p(\hat{\nabla},\hat{\mu}),$ the averages
$$A_{\mathbf{n}}(\widehat{f})=\frac{1}{|\nbb|}
\sum\limits_{\mathbf{k}=1}^{\mathbf{n}}\alpha(\mathbf{k})\mathbf{T}^{\mathbf{k}}(\widehat{f})
$$
$(o)$- converge to some $\widetilde{f}$ in
$L_p(\hat{\nabla},\hat{\mu})$.
\end{thm}

\begin{pf} Using Theorem \ref{T-o} we immediately obtain that
$${\mathbf{T}}^{{\mathbf{k}}}(\widehat{f})(\omega)
=(T_1^{k_1}\cdots
T_d^{k_d}\widehat{f})(\omega)=T_1^{k_1}(\omega)\cdots
T_d^{k_d}(\omega)f(\omega)={\mathbf{T}}^{\kb}(\omega)f(\omega)
 $$
 for any $\widehat{f}\in
 L_p(\hat{\nabla},\hat{\mu})$
and for  a.e.
 $\omega\in\Omega$.  Hence,
 $$A_{\mathbf{n}}(\widehat{f})(\omega)=\frac{1}{|\nbb|}
\sum\limits_{{\mathbf{k}}=1}^{{\mathbf{n}}}\alpha({\mathbf{k}}){\mathbf{T}}^{{\mathbf{k}}}(\omega)(f(\omega))=A_{\nbb}(f(\omega))
$$
for  a.e.
 $\omega\in\Omega$ and for any $\widehat{f}\in
 L_p(\hat{\nabla},\hat{\mu}).$

Since every $T_i$ is  linear  positive
$L_1(\hat{\nabla},\hat{\mu})$ contractions in
$L_p(\hat{\nabla},\hat{\mu})$ and
 $T_i\id\leq\id$ then by Theorem \ref{T-o} every
 $T_i(\omega)$ is positive $L_1(\nabla_\omega,\mu_\omega)$ contraction in $L_p(\nabla_\omega,\mu_\omega)$ and
 $T_i\id(\omega)\leq\id(\omega)$ for  a.e.
 $\omega\in\Omega$.  This means that $T_i(\omega)$ is
 $L_1(\nabla_\omega,\mu_\omega)-L^\infty(\nabla_\omega,\mu_\omega)$-
 contraction.  Then by Theorem 1.2 \cite{JO} the averages $A_\nbb(f(\omega))$
 $(o)$-converge to some  $g(\omega)\in
 L_p(\nabla_\omega,\mu_\omega).$ According Corollary \ref{2.5} we find
 $A_{\mathbf{n}}(\widehat{f})=\widehat{A_{\mathbf{n}}(f(\omega))}
 \stackrel{(o)}\rightarrow\widetilde{f}=\widehat{g(\omega)}$ in $L_0(\hat{\nabla},\hat{\mu}).$

From
\begin{eqnarray*}
|A_\nbb(f(\omega))|&=&|\frac{1}{|\nbb|}
\sum\limits_{{\mathbf{k}}=1}^{{\mathbf{n}}}\alpha(\kb){\mathbf{T}}^{\kb}(\omega)(f(\omega))|\\[2mm]
&\leq& \frac{1}{|\nbb|}
\sum\limits_{{\mathbf{k}}=1}^{{\mathbf{n}}}|\alpha({\mathbf{k}})|{\mathbf{T}}^{{\mathbf{k}}}(\omega)(|f(\omega)|)\\[2mm]
&\leq & \frac{b}{|\nbb|}
\sum\limits_{{\mathbf{k}}=1}^{{\mathbf{n}}}{\mathbf{T}}^{{\mathbf{k}}}(\omega)(|f(\omega)|)
\end{eqnarray*}
and
$$\sup\limits_{{\mathbf{n}}}\frac{1}{|\nbb|}\sum\limits_{{\mathbf{k}}=1}^{{\mathbf{n}}}{\mathbf{T}}^{{\mathbf{k}}}(\omega)(|f(\omega)|)\in
L_p(\nabla_\omega,\mu_\omega),
$$ we find
$\sup\limits_{{\mathbf{n}}}|A_{\mathbf{n}}(f(\omega))|\in
L_p(\nabla_\omega,\mu_\omega),$ where
$b=\sup\limits_{{\mathbf{k}}}|\alpha(\kb)|$ for  a.e.
$\omega\in\Omega.$

According Theorem \ref{2.3} one gets
$\sup\limits_{{\mathbf{n}}}|A_{\mathbf{n}}(\widehat{f})|\in
L_0(\hat{\nabla},\hat{\mu}).$

From Theorem 1.2 \cite{Kr}(p.196) it follows that
$$\bigg\|\sup\limits_{\nbb}|A_\nbb(f(\omega))|\bigg\|_p\leq
b\bigg\|\sup\limits_{\nbb}\frac{1}{|\nbb|}
\sum\limits_{\kb=1}^{\nbb}{\mathbf{T}}^{{\mathbf{k}}}(\omega)(|f(\omega)|)\bigg\|_{L_p(\nabla_\omega,\mu_\omega)}\leq
b \cdot q^d \|f(\omega)\|_{L_p(\nabla_\omega,\mu_\omega)}$$ for
a.e. $\omega\in\Omega.$

This means that
$$\|\sup\limits_{{\mathbf{n}}}|A_{\mathbf{n}}(\widehat{f})||_{L_p(\hat{\nabla},\hat{\mu})}\leq b \cdot q^d
\|\widehat{f}\|_{L_p(\hat{\nabla},\hat{\mu})}$$ and
$A_{\mathbf{n}}(\widehat{f})$ is bounded in
$L_p(\hat{\nabla},\hat{\mu}).$ Hence
$A_{\mathbf{n}}(\widehat{f})\stackrel{(o)}\rightarrow\widetilde{f}$
in $L_p(\hat{\nabla},\hat{\mu}).$
\end{pf}

\section*{Acknowledgement}
The authors are grateful to Professor Vladimir Chilin for his
valuable comments and remarks on improving the paper. The secon
named author (F.M.) acknowledges the MOHE Grant FRGS11-022-0170. He
also thanks the Junior Associate scheme of the Abdus Salam
International Centre for Theoretical Physics, Trieste, Italy.

\end{document}